\documentclass[11pt]{article}
\usepackage{amsmath}
\usepackage{amssymb}
\usepackage{latexsym}

\usepackage{amsthm}
\theoremstyle{plain}
\newtheorem{thm}{Theorem}[section]   
\newtheorem{prop}[thm]{Proposition}

\newtheorem{lem}[thm]{Lemma}

\theoremstyle{definition}
\newtheorem{df}[thm]{Definition}

\newtheorem{exs}[thm]{Examples}
\newcommand{\pf}{\noindent{\sc Proof.}\ }

\newcommand{\pfend}%
{\hspace*{\fill}\lower3pt\hbox{$\Box$}\medskip} 

\usepackage{subfig}
\usepackage{caption}
\newcommand{\R}{\mathbb{R}}

\renewcommand{\phi}{\varphi}
\newcommand{\Ga}{\Gamma}
\newcommand{\Om}{\Omega}

\newcommand{\add}[1]{\mathbin{\lower 5pt%
    \hbox{${\stackrel{\textstyle +}{\scriptscriptstyle #1}}$}}}

\usepackage{xspace}

\newcommand{\VBgpd}{$\mathcal{V}\!\mathcal{B}$--group\-oid\xspace} 
\newcommand{\VBgpds}{$\mathcal{V}\!\mathcal{B}$--group\-oids\xspace} 
\newcommand{\LAgpd}{$\mathcal{L}\!\mathcal{A}$--group\-oid\xspace}
\newcommand{\LAgpds}{$\mathcal{L}\!\mathcal{A}$--group\-oids\xspace}
\newcommand{\LAvb}{$\mathcal{L}\!\mathcal{A}$--vector--bundle\xspace}

\newcommand{\VBalgd}{$\mathcal{V}\!\mathcal{B}$--algebr\-oid\xspace}

\renewcommand{\leq}{\leqslant}

\newcommand{\st}{\ \vert\ }

\renewcommand{\Bar}[1]{\overline{#1}}

\newcommand\co{\colon\thinspace}

\newcommand{\Dd}{D}

\newcommand{\gog}{\mathfrak{g}}

\newcommand{\sol}{\bullet}

\usepackage{mathrsfs}

\newcommand{\cinf}[1]{C^{\infty}(#1)}
\newcommand{\sfn}[1]{C^{\infty}(#1)}
\newcommand{\vf}[1]{\mscr{X}(#1)}

\newcommand{\act}{\mathbin{\hbox{$<\kern-.4em\mapstochar\kern.4em$}}}
\newcommand{\ract}{\mathbin{\hbox{$\mapstochar\kern-.3em>$}}}

\usepackage{mathrsfs}
\newcommand{\mscr}[1]{\mathscr{#1}}

\newcommand{\sg}{\mathscr{G}}
\newcommand{\scrh}{\mathscr{H}}
\newcommand{\gpd}{\rightrightarrows}

\newcommand{\ld}{\mathfrak{L}}             % Lie derivative

\newcommand{\Ri}[1]{\overrightarrow{#1}}

\newcommand{\tilalpha}{\widetilde\alpha}   
\newcommand{\tilbeta}{\skew6\widetilde\beta}
\newcommand{\tilc}{\widetilde{c}}   
\newcommand{\hato}{\widehat{0}}
\newcommand{\tilo}{\widetilde{0}}   
\newcommand{\tilq}{\widetilde{q}}

\newcommand{\Rarr}{\overrightarrow}

\newcommand{\sdp}{\ltimes}

\newcommand{\D}{\mathop{\raise0.1ex\hbox{$\mathfrak{D}$}}} 

\newcommand{\End}{\mathop{\rm End}}

\newcommand{\Ad}{{\rm Ad}} % Do NOT use Declaremathop

\usepackage[all]{xy}

\newdir{ >}{{}*!/-5pt/@{>}} % exercise 14 in short doc

\newcommand{\ses}[3]{\SelectTips{cm}{12}%
\xymatrix@=8mm{#1\ar@{ >->}[r] & #2 \ar@{->>}[r] & #3}}

\newcommand{\extt}[2]{\mathsf{\Lambda}^{#1}(#2)}

\newcommand{\tilone}{\widetilde 1}

\usepackage{array}

\usepackage{url}

\begin{document}

\title{\textbf{From symplectic groupoids to double structures}}

\author{Kirill Mackenzie\\[2mm]
School of Mathematics and Statistics,\\
University of Sheffield,\\
Sheffield S3 7RH, UK\\
\url{K.Mackenzie@sheffield.ac.uk}}

\maketitle

\tableofcontents

\newpage

\section*{Foreword}
\addcontentsline{toc}{section}{Foreword}

These notes are an introduction to symplectic groupoids and the double structures
associated with them. The treatment is intended to lie about midway between the 
original account of Coste, Dazord and Weinstein \cite{CosteDW}, which relied on 
effective use of the symplectic structures, and the account in my book 
\cite{Mackenzie:GT}, which showed, on the level of Poisson groupoids,  
that the basic results of the theory follow from `categorical' compatibility 
conditions between the associated Lie algebroid and Lie groupoid structures. 
(See \S\ref{sect:da} for more details.) 

The reader needs to know only the most basic ideas of symplectic geometry, Lie
groups, and vector bundles. Conventions are recalled in \S\ref{sect:back}. 

In particular, no familiarity with double structures is assumed. Instead the
cotangent groupoid --- perhaps one of the hardest ideas to assimilate for someone 
new to this theory --- is introduced as a consequence of the isomorphism between  
the tangent and cotangent bundles. That then leads to the general concept of a
\VBgpd and its duality. 

These notes are based on lectures given at the School on 
\emph{Geometric, Algebraic and Topological Methods for Quantum
Field Theory} in Villa de Leyva, Colombia, in July 2013. 
I am very glad to have had the opportunity to give these
lectures, and I thank the organizers most heartily for involving me in a 
School of such vitality and openness. 
I want to particularly thank Alexander Cardona for his 
generous hospitality and for looking after me so well.  

These notes aim to introduce the reader to certain important, and 
relatively new, ideas quickly; accordingly they omit much standard 
material. All the main results here are known, but the approach has 
some new features. Some references which provide alternative 
treatments and more detail are given at the end. 

I am very grateful to Yvette Kosmann--Schwarzbach for her careful
reading and valuable comments. 

\section{Background: Poisson structures, Lie algebroids, Lie groupoids} 
\label{sect:back}

\begin{df}
A \emph{Poisson structure} on a manifold $P$ is a bracket 
of smooth functions 
$$
\{\ ,\ \}\co \sfn{P}\times \sfn{P}\to \sfn{P}
$$ 
with respect to which $\sfn{P}$ is an $\R$\negthinspace--Lie 
algebra, and such that for all $f_1, f_2, f_3\in \sfn{P},$
\begin{equation}                           
\label{eq:Pbrack}
\{f_1, f_2f_3\} = f_2\{f_1,f_3\} + \{f_1,f_2\}f_3.
\end{equation}
\end{df}

The bracket $\{f_1, f_2\}$ depends only on $df_1$ and $df_2$ and we define a
2-vector field $\pi$ on $P$ by $\pi(w_1\,df_1, w_2\,df_2) = w_1w_2\{f_1, f_2\}$.
This is the \emph{Poisson tensor} and in practice we use $\pi$ to denote 
a Poisson structure. 

For a symplectic manifold $(M,\omega)$, there is an \emph{associated 
Poisson structure} on $M$ defined by $\{f_1,f_2\} = \omega((df_1)^\sharp,(df_2)^\sharp)$,
where $\phi\mapsto\phi^\sharp,$ $T^*M\to TM$ is the inverse of the map 
$\omega^\flat\co TM\to T^*M,$
$\langle\omega^\flat(X),Y\rangle = -\omega(X,Y)$. 
% THIS SIGN INTRODUCED FOR COMPATIBILITY WITH LATER

Returning to a Poisson manifold $P$, we define $\pi^\#\co T^*P\to TP$, the 
\emph{Poisson anchor}, by 
\begin{equation}
\langle \phi, \pi^\#(\psi)\rangle = \pi(\psi, \phi).   
\end{equation}
If $\pi^\#$ is an isomorphism of vector bundles, then 
$\omega(X,Y) = \langle (\pi^\#)^{-1}(Y), X\rangle$ % CHECK !! 
is a symplectic structure on $P$, and the associated Poisson structure is the
given one. 

Define a bracket of 1-forms by
\begin{equation}
\label{eq:ldf}
[\phi, \psi] = \ld_{\phi^\#}(\psi)  -  \ld_{\psi^\#}(\phi) 
  - d(\pi(\phi, \psi)),
\end{equation}
where again $\phi^\# = \pi^\#(\phi)$. The bracket of 1-forms makes $T^*P$ a Lie 
algebroid with anchor $\pi^\#$. 

\begin{df}
A \emph{Lie algebroid} is a vector bundle $A$ on base $M$ together with a bracket
$[\ ,\ ]\co \Ga A\times\Ga A\to \Ga A$, with respect to which $\Ga A$ is an
$\R$-Lie algebra, and a vector bundle map $a\co A\to TM$, called the \emph{anchor 
of $A$}, such that 
\begin{equation}
[X,fY] = f[X,Y] + a(X)(f)Y,\qquad a[X,Y] = [aX, aY],    
\end{equation}
for $X,Y\in\Ga A$ and $f\in\sfn{M}$. 

If $A'$ is also a Lie algebroid on $M$ then a \emph{base-preserving morphism of Lie 
algebroids} is a morphism of vector bundles $\phi\co A\to A'$ such that $a'\circ\phi = a$
and $\phi([X,Y]) = [\phi(X), \phi(Y)]$ for all $X, Y\in\Ga A$. 
\end{df}

For any vector bundle $A$ on $M$, sections $X\in\Ga A$ induce (fibrewise) linear maps 
$\ell_X\co A^*\to\R$ by $\ell_X(\phi) = \langle\phi, X(q_*(\phi))\rangle$ where $q_*$ is
the bundle projection for $A^*$. Every fibrewise linear map $A^*\to\R$ is $\ell_X$ for
a unique $X\in\Ga A$. 

When $A$ is a Lie algebroid, there is a unique Poisson structure on the manifold $A^*$, 
called the \emph{Lie-Poisson structure corresponding to $A$}, 
such that
\begin{equation}
\label{eq:L-P}
\{\ell_X, \ell_Y\} = \ell_{[X,Y]},\qquad   
\{\ell_X, f\circ q_*\} = a(X)(f)\circ q_*,\qquad
\{f_1\circ q_*, f_2\circ q_*\} = 0, 
\end{equation}
for $X,Y\in\Ga A$, $f\in\sfn{M}$. It is not clear that this defines the bracket 
for all pairs of smooth functions. To bypass this, consider the Poisson tensor 
for $A^*$, and note that every 1-form is the sum of $d\ell_X$ for some $X\in\Ga A$ 
and $q_*^*df$ for some $f\in\sfn{M}$. Conversely, if $E$ is a vector bundle on
$M$ with a Poisson structure for which: (i) the bracket of fibrewise linear functions
is fibrewise linear, (ii) the bracket of a fibrewise linear function with a pullback 
function is a pullback function, and (iii) the bracket of two pullback functions is zero,
then there is a unique Lie algebroid structure on $E^*$ for which this is the 
Lie-Poisson structure. 

For a Lie algebra $\gog$, the second and third definitions in (\ref{eq:L-P})
are vacuous. However the first can be extended to any pair of smooth functions. 
Let $f\co\gog^*\to\R$ be a smooth function on $\gog^*$. Using an elementary notion
of derivative, $D(f)(\phi)$, for each $\phi\in\gog^*$, is a linear map $\gog^*\to\R$
(the `directional derivative') and can thus be identified with an element of $\gog$. 
Since a Poisson bracket depends only on the derivatives of its arguments, we
have
\begin{equation}
\{f_1, f_2\}(\phi) = \langle\phi,[D(f_1)(\phi), D(f_2)(\phi)]\rangle.   
\end{equation}

A map $\mu\co P\to Q$ of Poisson manifolds is a \emph{Poisson map} if 
$\{f_1\circ \mu, f_2\circ \mu\}_P = \{f_1, f_2\}_Q\circ \mu$ for all 
$f_1, f_2\in\cinf{Q}$. We will need Poisson maps in \S\ref{sect:pgpds}.

\subsection*{Lie groupoids}

\begin{df}
A \emph{Lie groupoid} consists of a manifold $M$ of \emph{points}, 
or \emph{objects}, and a manifold $\sg$ of \emph{arrows}, together with surjective 
submersions $\alpha, \beta\co\sg\to M$, called respectively the \emph{source} and 
\emph{target} maps, and an injective immersion $1\co M\to\sg$, written $m\mapsto 1_m$, 
together with a smooth \emph{multiplication} or \emph{composition}, $(h,g)\mapsto hg$, 
$\sg\times_M\sg\to\sg$, defined on 
$$
\sg\times_M\sg = \{(h,g)\in\sg\times\sg\st \alpha(h) = \beta(g)\}. 
$$
These maps are subject to modified forms of the group axioms: for all
$h,g,f\in\sg$ such that $\alpha(h) = \beta(g)$ and $\alpha(g) = \beta(f)$, 
we have $h(gf) = (hg)f$; for each $g\in\sg$ we have $1_ng = g$ and $g1_m = g$ 
where $n = \beta(g)$, $m = \alpha(g)$; finally, for each $g\in\sg$ there is a
$g^{-1}\in\sg$ such that $g^{-1}g = 1_m$ and $gg^{-1} = 1_n$. 
We often write $\sg\gpd M$ to indicate a Lie groupoid. 

A \emph{morphism of Lie groupoids} consists of two smooth maps $F\co \sg\to\sg'$
and $f\co M\to M'$ which commute with the sources and targets, for which
$F(1_m) = 1'_{f(m)}$ for $m\in M$, and for which $F(hg) = F(h)F(g)$ when
$\alpha(h) = \beta(g)$. If $M' = M$ and $f$ is the identity map, we say that $F$ 
is \emph{over $M$}.
\end{df}

In these lectures we are concerned with symplectic groupoids, rather than the
general theory. We include here only four important examples. 

\begin{exs}
{\bf (i)} For any manifold $M$ the Cartesian square $M\times M$ is a Lie groupoid with 
$\alpha(n,m) = m$, $\beta(n,m) = n$, $1_m = (m,m)$ and multiplication
$(p,n)(n,m) = (p,m)$. This is the \emph{pair groupoid on $M$}.   

\smallskip

{\bf (ii)} For a connected manifold $M$ the \emph{fundamental groupoid} $\Pi(M)$ is a
Lie groupoid on $M$. The arrows are homotopy classes, with endpoints fixed, 
of piecewise smooth paths in $M$, the source and target maps give the start-point
and end-point, and the groupoid multiplication is concatentation. The identity
elements are the classes of constant paths. The map $\chi\co \Pi(M)\to M\times M$
which sends the class of a path $p(t)$ to $(p(1), p(0))$ is a morphism of Lie
groupoids over $M$. 

\smallskip

{\bf (iii)} Let $E$ be a vector bundle on $M$. Denote by $\Phi(E)$ the set of all linear 
isomorphisms $\xi\co E_m\to E_n$ between fibres of $E$. The local triviality
of $E$ can be used to give $\Phi(E)$ a manifold structure. Define source and
target $\alpha(\xi) = m$ and $\beta(\xi) = n$, for $\xi$ as above, and write
$1_m$ for the identity map $E_m\to E_m$. Then the standard composition of
maps makes $\Phi(E)$ a Lie groupoid on base $M$, the \emph{frame groupoid of
$E$}. 

\smallskip

{\bf (iv)} Let $G$ be a Lie group and let $G\times M\to M$ be a smooth action of $G$ on
a manifold $M$. Write $\sg = G\times M$, the product manifold. Define
$\alpha\co \sg\to M$ by $\alpha(g,m) = m$ and $\beta\co\sg\to M$ by
$\beta(g,m) = gm$. For $m\in M$ write $1_m = (1,m)$ and define
$$
(h,n)(g,m) = (hg,m)
$$
when $n = gm$. Then $\sg$ is a Lie groupoid on base $M$, the \emph{action
groupoid} corresponding to the given action. We write $G\act M$ for $\sg$;
many authors write $G\sdp M$. 
\end{exs}

We now describe the Lie algebroid of a Lie groupoid. The construction follows 
the same general procedure as the construction of the Lie algebra of a Lie group,
though we define the bracket using right-invariant vector fields rather than 
the more usual left-invariant vector fields. 

For a Lie group $G$ the Lie algebra $\gog$ may be defined as $T_1(G)$, the tangent
space at the identity element, with the bracket on $\gog$ defined by identifying 
elements of $\gog$ with right-invariant vector fields. 

For a Lie groupoid $\sg\gpd M$ there is a regular submanifold of identity elements 
$\{1_m\st m\in M\}\subseteq\sg$, 
which we identify with $M$. Write $T_M\sg$ for the restriction of the tangent bundle
$T\sg$ to $M$. The fibres of $T_M\sg$ are $T_{1_m}\sg$ for $m\in M$. 

Because the multiplication in $\sg$ is only defined on $\sg\times_M\sg$, right-translations
are only defined on the fibres of the source map. For $g\in\sg$ we have
$$
R_g\co \alpha^{-1}(\beta g) \to \alpha^{-1}(\alpha g),\quad h\mapsto hg. 
$$
In order to define right-invariance for vector fields on $\sg$, we therefore consider only
vector fields which are vertical with respect to the source map. Thus:

\begin{df}
A vector field $\xi$ on $\sg$ is \emph{right-invariant} if $T(\alpha)(\xi(g)) = 0$ for all
$g\in\sg$, and $T(R_g)(\xi(h)) = \xi(hg)$ for all $(h,g)\in \sg\times_M\sg$.   
\end{df}

Note that for a smooth map $f\co M\to N$ we usually write $T(f)$ for the tangent 
map $TM\to TN$, reserving $df$ for when $f$ is real-valued. 

The second condition is equivalent to $\xi(g) = T(R_g)(\xi(1_{\beta g}))$ for all 
$g\in\sg$. If the analogy with the Lie algebra of a Lie group is to hold good, it
should be sufficient to consider the restrictions of right-invariant vector fields
to the manifold of identity elements. We therefore define
$$
A\sg = \ker (T(\alpha)\co T\sg\to TM) \cap T_M\sg. 
$$
That is, $A\sg$ is a vector bundle on $M$, and the fibre over $m\in M$ is
$T_{1_m}(\alpha^{-1}(m))$. 

If $X$ is a section of $A\sg$ we can define a right-invariant vector field 
$\Rarr{X}$ on $\sg$ by 
\newline
$\Rarr{X}(g) = T(R_g)(X(\beta g))$. Conversely, every
right-invariant vector field  is $\Rarr{X}$ for some $X\in\Ga A\sg$. 

We next need to show that the bracket of right-invariant vector fields is a 
right-invariant vector field. This follows, as in the case of Lie groups, from
the fact that the bracket of projectable vector fields is projectable: note 
that a vector field $\xi$ on $\sg$ is right-invariant if and only if
it projects to zero under $\alpha$ and for each $g\in\sg$, the restriction
of $\xi$ to $\alpha^{-1}(\beta g)$ projects under $R_g$ to the restriction 
of $\xi$ to $\alpha^{-1}(\alpha g)$. We therefore have a bracket on $\Ga A\sg$
such that, for all $X, Y\in\Ga A\sg$, 
$$
\Rarr{[X,Y]} = [\Rarr{X}, \Rarr{Y}].
$$
Lastly, consider $[X,fY]$ where $f\in\sfn{M}$. Clearly 
$\Rarr{fY} = (f\circ\beta)\Rarr{Y}$. Using the Leibniz identity for vector fields on
$\sg$, we therefore get that
$$
[X, fY] = f[X,Y] + a(X)(f)Y,
$$
where $a\co A\sg\to TM$ is the restriction of $T(\beta)\co T\sg\to TM$, called
the \emph{anchor of $A\sg$}. 
This completes the construction of the \emph{Lie algebroid of the Lie groupoid $\sg$.}

\begin{exs}
{\bf (i)} Consider $\sg = M\times M$. The kernel of $T(\alpha)$ is $\beta^!(TM)$, the 
pullback of $TM$ to $M\times M$ across $\beta$. Restricting this to the diagonal in
$M\times M$ we get $A\sg = TM$ as vector bundles. The right-invariant vector field
corresponding to $X\in\Ga TM$ is $\Rarr{X}(n,m) = X(n)\oplus 0_m$ and so the bracket
on $A\sg$ is the standard bracket of vector fields on $M$ and the anchor is the
identity. 

\smallskip

{\bf (ii)} It is easy to see that a base-preserving morphism of Lie groupoids 
$F\co\sg \to \sg'$ induces a map $A(F)\co A\sg \to A\sg'$ which is a morphism of
Lie algebroids. In the case of $\chi\co \Pi(M)\to M\times M$ the kernel (that is,
the set of elements of $\Pi(M)$ which are mapped to identity elements of $M\times M$), 
is the union of the fundamental groups of $M$ and is fibre-wise discrete. It follows
that $A(\chi)\co A\Pi(M)\to TM$ is a diffeomorphism, usually % thanks YKS
identified with the identity. 

\smallskip

{\bf (iii)} The Lie algebroid $A\Phi(E)$, as a vector bundle, is an extension of
$TM$ by $\End(E)$. Its sections are those linear differential operators $D$, of order
$\leq 1$, for which there exists a vector field $X$ on $M$ such that
$$
D(f \mu) = fD(\mu) + X(f)\mu,
$$
for all $f\in\sfn{M}$ and $\mu\in\Ga E$. These may be called \emph{derivations}
on $E$, and $A\Phi(E)$ denoted by $\D(E)$. The bracket on $\D(E)$ is the standard
commutator bracket, and the anchor maps each $D$ to its $X$. 

We will not prove this, but it is easy to believe: the sections of $\D(E)$ are all
covariant derivatives $\nabla_X$ for all connections $\nabla$ and all vector fields
$X$ on $M$, together with their differences, and $\Phi(E)$ consists of all parallel
translations, for all connections in $E$ and all paths in $M$. 

\smallskip

{\bf (iv)} Let $G\times M\to M$ be a smooth action of a Lie group $G$ on a manifold 
$M$ and let $X\mapsto X_M$ be the infinitesimal action $\gog\to \vf{M}$. 

The Lie algebroid of $G\act M\gpd M$ is the trivial vector bundle
$M\times\gog$ with anchor $M\times\gog\to TM$, $(m,X)\mapsto X_M(m)$ and bracket
$$
[V,W] = L_{V_M}(W) - L_{W_M}(V) + [V,W]^{\text{pt}}.
$$
Here $V$ and $W$ are $\gog$-valued maps on $M$; since $\gog$ is a vector space we
can take Lie derivatives of $\gog$-valued maps. The final term is the `pointwise
bracket' of $V$ with $W$ as maps into $\gog$. 
\end{exs}

\section{Symplectic groupoids}

The following definition is due to Weinstein \cite{Weinstein:1987}. 
Independent definitions were given by Karas\"ev \cite{Karasev:1986} 
and Zakrzewski \cite{Zakrzewski:1990-I,Zakrzewski:1990-II}. 

\begin{df}
\label{df:symgpd}
A \emph{symplectic groupoid} is a Lie groupoid $\Sigma\gpd P$ together with a symplectic
structure $\omega$ on $\Sigma$ such that the graph of multiplication
\begin{equation}
\label{eq:graph}
Gr = \{(hg,h,g)\st \alpha h = \beta g\}
\end{equation}
is Lagrangian in $\Bar{\Sigma}\times\Sigma\times\Sigma$. The bar denotes reversal of the
symplectic structure. 
\end{df}

The following are the fundamental consequences:
\begin{thm}
\renewcommand{\theenumi}{\rm(\roman{enumi})}
\label{thm:sgpd}
\begin{enumerate}
\item $1_P$ is Lagrangian in $\Sigma$.
\item Inversion is an antisymplectomorphism.
\item There is a unique Poisson structure on $P$ such that $\beta$ is a Poisson map.
\item There is a canonical isomorphism of Lie algebroids $A\Sigma\cong T^*P$. 
\end{enumerate}
\end{thm}

The proof occupies the rest of the section. First the fundamental example. 

\subsection*{The cotangent groupoid of a Lie group}

Consider any Lie group $G$. The tangent bundle $TG$ has a Lie group structure with
multiplication $T(\kappa)\co TG\times TG\to TG$ where $\kappa\co G\times G\to G$ is the
multiplication in $G$. Since the domain is a product manifold, we can differentiate
in each variable separately and we get
$T(\kappa)(Y,X) = T(L_h)(X) + T(R_g)(Y)$ where $Y\in T_hG$ and $X\in T_gG$. 
In practice we write $Y\sol X  = T(\kappa)(Y,X)$. With this structure, $TG$ is
the \emph{tangent group of} $G$, 

It is possible to give the cotangent bundle $T^*G$ a group structure in a similar
way. However what we are about to define is distinct from the group structure. 

Define a groupoid structure on $T^*G$ with base $\gog^*$ by
\begin{equation}
\label{eq:T*G}
\begin{split}
\tilbeta(\phi_g) = \phi\circ T_1(R_g), \qquad
\tilalpha(\phi_g) = \phi\circ T_1(L_g), \\
\psi_h\cdot\phi_g = \psi\circ T(R_{g^{-1}}) = \phi\circ T(L_{h^{-1}}).\phantom{XX}
% Corrected April 15, 2014
\end{split}
\end{equation}

This structure is isomorphic to an action groupoid. Let $G$ act on $\gog^*$ by the
coadjoint action $g\theta := \theta\circ \Ad_{g^{-1}}$ and form $G\act\gog^*\gpd\gog^*$. 
Then $G\act\gog^*\to T^*G$, $(g,\theta)\mapsto \theta\circ T(L_{g^{-1}})$ is a groupoid
isomorphism. 

Whereas the group $TG$ is isomorphic to the semi-direct product group $G\sdp\gog$
defined by the adjoint representation, $T^*G$ is isomorphic to the action groupoid 
$G\act\gog^*$. It is important to distinguish the two constructions.  

The symplectic structure on $T^*G$ is the standard symplectic structure for a cotangent 
bundle, $\omega = d\lambda$ where $\lambda$ is the 1-form on $G$
defined on $\xi\in T(T^*G)$ by 
$$
\langle\lambda,\xi\rangle = \langle\phi,T(c)(\xi)\rangle
$$
where $\phi\in T^*G$ is the base-point of $\xi$ and $c\co T^*G\to G$ is the projection.

Further, $\omega$ makes $T^*G\gpd\gog^*$ a symplectic groupoid and % that
the Poisson structure on $\gog^*$ which arises from Theorem \ref{thm:sgpd}(iii) is
the Lie--Poisson structure dual to $\gog$. 

\subsection*{Proof of Theorem \ref{thm:sgpd}}

To see what Definition \ref{df:symgpd} means, we use the tangent 
groupoid of $\Sigma$. This construction applies to any Lie groupoid $\sg\gpd M$. 
The \emph{tangent groupoid} of $\sg\gpd M$ is the groupoid $T\sg\gpd TM$
obtained by applying the tangent functor to the structure of $\sg\gpd M$. So the 
source of the tangent groupoid is $T(\alpha)$, the identities are
$T(1)(x)$ for $x\in TM$, and the multiplication $\eta\bullet \xi$ of two tangent vectors 
with $T(\alpha)(\eta) = T(\beta)(\xi)$ is $T(\kappa)(\eta,\xi)$ where 
$\kappa\co\sg\times_M\sg\to\sg$ is the multiplication in $\sg$. Because the tangent 
functor preserves diagrams, 
the groupoid axioms follow immediately. We write $\xi^{-1}$ for the inverse of $\xi\in T\sg$, 
and $\tilo_g$ for the zero element of $T_g\sg$. 

If $\sg = G$ is a group, then, as we saw above, $TG$ is a group with 
$\eta\bullet \xi = T(R_g)(\eta) + T(L_h)(\xi)$ for $\eta\in T_h(G)$ 
and $\xi\in T_g(G)$. However no Lie group has a symplectic groupoid structure.
For general Lie groupoids, the following result is very necessary. 

\begin{prop}
\label{prop:Tk}
For elements $\xi_i\in T\sg$, $i = 1,\dots,4$ with $\xi_1, \xi_2\in T_g\sg$
and $\xi_3, \xi_4\in T_h\sg$, and with $T(\alpha)(\xi_3) = T(\beta)(\xi_1)$ and
$T(\alpha)(\xi_4) = T(\beta)(\xi_2)$, we have
\begin{equation}
\label{eq:Tk}
(\xi_4+\xi_3)\sol(\xi_2+\xi_1) = \xi_4\sol \xi_2 + \xi_3\sol \xi_1  
\end{equation}
\end{prop}

\pf
This is the statement that $T(\kappa)$ is linear as a map of vector bundles
$T\sg\times_{TM} T\sg\to T\sg$. 

To see this, observe that $(\xi_3, \xi_1)$ and $(\xi_4, \xi_2)$ are elements of
$T_h\sg\times T_g\sg$ and can therefore be added to give 
$(\xi_4+\xi_3,\xi_2+\xi_1)\in T_h\sg\times T_g\sg$.

The linearity then gives that
$$
T(\kappa)(\xi_4+\xi_3,\xi_2+\xi_1) = T(\kappa)(\xi_4, \xi_2) + T(\kappa)(\xi_3, \xi_1)
$$
which is (\ref{eq:Tk}). 
\pfend

The source and target conditions ensure that each product and sum is defined and 
is in the domain of $T(\kappa)$. 

Consider now a symplectic groupoid $\Sigma\gpd P$. 
The tangent to the graph (\ref{eq:graph}) of the multiplication is the graph of the 
tangent multiplication~:
\begin{equation}
\label{eq:tgraph}
T(Gr) = \{(\eta\sol \xi,\eta,\xi)\st T(\alpha)(\eta) = T(\beta)(\xi)\}.
\end{equation}
Recall that a submanifold $L$ of a symplectic manifold $(S, \sigma)$ is \emph{Lagrangian} 
if it is isotropic: $\sigma(Y,X) = 0$ for all $X,Y\in TL$, and $\dim L = \frac{1}{2}\dim S$.

Applying the isotropy condition to the graph of the multiplication we have that
\begin{equation}
\label{eq:isot}
-\omega(\eta_1\sol \xi_1, \eta_2\sol \xi_2) + \omega(\eta_1, \eta_2) 
+ \omega(\xi_1, \xi_2) = 0
\end{equation}
for all $\xi_1, \xi_2, \eta_1, \eta_2\in T\Sigma$ for which the multiplications are defined. 

Set $\xi_i = \eta_i = T(1)(x_i)$ in (\ref{eq:isot}), where $x_i\in TP$. Then (\ref{eq:isot}) 
becomes
$$
-\omega(T(1)(x_1), T(1)(x_2)) + \omega(T(1)(x_1), T(1)(x_2)) + 
\omega(T(1)(x_1), T(1)(x_2)) = 0
$$
so $\omega(T(1)(x_1), T(1)(x_2)) = 0$ for all $x_1, x_2\in TP$. 
This proves that $1_P$ is isotropic in $\Sigma$. Next, 
$$
\dim (Gr) = \dim(\Sigma\times_P\Sigma) = 2\dim\Sigma - \dim P
$$
and this is to be $\frac{1}{2}\dim (\Sigma\times\Sigma\times\Sigma)$ 
so $\dim P = \frac{1}{2}\dim \Sigma$. 
This completes the proof that $1_P$ is Lagrangian in $\Sigma$. 

Next consider any $\xi_1, \xi_2$ and write $x_i = T(\alpha)(\xi_i)$. Then with 
$\eta_i = \xi_i^{-1}$ in (\ref{eq:isot}) we have
$$
-\omega(T(1)(x_1), T(1)(x_2)) + \omega(\xi_1^{-1}, \xi_2^{-1}) + 
\omega(\xi_1, \xi_2) = 0
$$
So 
\begin{equation}
\label{eq:inv}
\omega(\xi_1^{-1}, \xi_2^{-1}) = -\omega(\xi_1, \xi_2), 
\end{equation}
which shows that inversion is antisymplectic. 

The last two parts of Theorem \ref{thm:sgpd} require more work. 

Recall the map $\omega^\flat\co T\Sigma\to T^*\Sigma$ defined by 
$\omega^\flat(\xi)(\eta) = -\omega(\xi,\eta).$
What is $\omega^\flat(\eta\sol \xi)$, where $\eta\in T_h\Sigma$, $\xi\in T_g\Sigma$
and $T(\alpha)(\eta) = T(\beta)(\xi)$~? 
It is an element of $T^*_{hg}\Sigma$ so we pair it with an arbitrary element $\zeta$
of $T_{hg}\Sigma$.

We can write $\zeta = \zeta_2\sol \zeta_1$ where $\zeta_2\in T_h\Sigma$ and 
$\zeta_1\in T_g\Sigma$. 
To see this, take any $\zeta_1\in T_g\Sigma$ with $T(\alpha)(\zeta_1) = T(\alpha)(\zeta)$ 
and define $\zeta_2 = \zeta\sol \zeta_1^{-1}$. This decomposition is not unique, of course. 

Then, by (\ref{eq:isot}), 
$$
\omega(\eta\sol \xi, \zeta_2\sol \zeta_1) = \omega(\eta, \zeta_2) + \omega(\xi, \zeta_1).
$$
Is the RHS well-defined~? Suppose we also have $\zeta = \zeta_4\sol \zeta_3$ where 
$\zeta_4\in T_h\Sigma$ 
and $\zeta_3\in T_g\Sigma$. Then $\zeta_3\sol \zeta_1^{-1} = 
\zeta_4^{-1}\sol \zeta_2\in T_{1_m}\Sigma$ where
$m = \alpha h = \beta g$. Call this $\nu$. So $\zeta_3 = \nu\sol \zeta_1$ and 
$\zeta_4 = \zeta_2\sol \nu^{-1}$. Now 
$$
\omega(\eta\sol \xi, \zeta_4\sol \zeta_3) = 
\omega(\eta, \zeta_2\sol \nu^{-1}) + \omega(\xi, \nu\sol \zeta_1).
$$
Consider the second term on the RHS first. Inserting an identity element and
using (\ref{eq:isot}), we have 
$$
\omega(\xi, \nu\sol \zeta_1) = \omega(T(1)T(\beta)(\xi)\sol \xi, \nu\sol \zeta_1) 
= \omega(T(1)T(\beta)(\xi), \nu) + \omega(\xi, \zeta_1). 
$$
Now the first term:
$$
\omega(\eta, \zeta_2\sol \nu^{-1}) = 
\omega(\eta\sol T(1)(T(\alpha)(\eta), \zeta_2\sol \nu^{-1}) = 
\omega(\eta,\zeta_2) + \omega(T(1)(T(\alpha)(\eta)), \nu^{-1}).
$$
Lastly, using (\ref{eq:inv}),  
$$
\omega(T(1)(T(\alpha)(\eta)), \nu^{-1}) = -\omega(T(1)(T(\alpha)(\eta)), \nu).
$$
So 
\begin{equation}
\label{eq:wd}
\omega(\eta\sol \xi, \zeta) = \omega(\eta, \zeta_2) + \omega(\xi, \zeta_1)
\end{equation}
is well-defined. 

\subsubsection*{The cotangent groupoid $T^*\Sigma$}

Now we define a groupoid structure on $T^*\Sigma$ so that (\ref{eq:wd}) becomes
\begin{equation}
\omega^\flat(\eta\sol \xi) = \omega^\flat(\eta)\sol \omega^\flat(\xi). 
\end{equation}
In fact, we have already done most of the work. 
The base of $T^*\Sigma$ will be $A^*\Sigma$ and the source and target maps 
$T^*\Sigma\to A^*\Sigma$ are
\begin{equation}
\label{eq:cotst}
\langle\tilbeta(\Phi), Y\rangle = \langle\Phi, T(R_g)(Y)\rangle,
% X, Y here are correct !!!
\quad
\langle\tilalpha(\Phi), X\rangle = \langle\Phi, T(L_g)(X - T(1)(aX))\rangle,
\end{equation}
where $\Phi\in T^*_g\Sigma$ and $Y\in A_{\beta g}\Sigma,\ X\in A_{\alpha g}\Sigma.$ 

For $\Phi\in T_g^*\Sigma$ and $\Psi\in T_h^*\Sigma$ define $\Psi\sol\Phi\in T^*_{hg}\Sigma$ by 
\begin{equation}
\label{eq:add}
\langle \Psi\sol\Phi, \eta\sol \xi\rangle = 
\langle \Psi, \eta\rangle + \langle \Phi, \xi\rangle. 
\end{equation}
The proof that the RHS of (\ref{eq:add}) is well-defined needs precisely the 
condition that the source and target match, and follows the same pattern as the proof 
that (\ref{eq:wd}) is well-defined. We leave the reader to check that when $\Sigma
= T^*G$ for a Lie group $G$, (\ref{eq:add}) reduces to (\ref{eq:T*G}). 

To define the identity $\tilone_\phi\in T^*_{1_m}\Sigma$ corresponding to 
$\phi\in A_m^*\Sigma,$ observe that every element $\xi$ of $T_{1_m}\Sigma$ can be written 
uniquely in the form $T(1)(x) + X$ where $x\in T_mM$ and 
$X\in A_m\Sigma.$ We can therefore define
\begin{equation}
\label{eq:cotid}
\langle\tilone_\phi, T(1)(x) + X\rangle = \langle\phi, X\rangle
\end{equation}
and it is straightforward to check that $\tilone_\phi$ is indeed an identity
for the multiplication and that this structure makes $T^*\Sigma\gpd A^*\Sigma$ a
groupoid with inverses
\begin{equation}
\label{eq:inverses}
\langle\Phi^{-1}, \xi^{-1}\rangle = - \langle\Phi, \xi\rangle.
\end{equation}

We now prove that $\omega^\flat\co T\Sigma\to T^*\Sigma$ is a morphism of groupoids, 
\begin{equation}
\label{eq:omegaflat}
{\xymatrix@=1.2cm{
T\Sigma \ar[r]^{\omega^\flat}
  \ar@<-0.5ex>[d]\ar@<0.5ex>[d]
& T^*\Sigma \ar@<-0.5ex>[d]\ar@<0.5ex>[d]\\
TP \ar[r]_w & A^*\Sigma\\
}}
\end{equation}
over some map $w$. To find $w$, 
take the composite of $T(1)$, followed by $\omega^\flat$, followed by $\tilbeta$
(or $\tilalpha$). This gives 
$$
w(x) = -\left.\iota_{T(1)(x)}\omega\right|_{A\Sigma}.
$$
We must now show that $\tilbeta\circ \omega^\flat = w\circ T(\beta)$. 
Take $\xi\in T_g\Sigma$. Then
$$
\langle\tilbeta(\omega^\flat(\xi)), Y\rangle =
\langle\omega^\flat(\xi), T(R_g)(Y)\rangle =
-\omega(\xi, T(R_g)(Y)).
$$
On the other hand, 
$$
\langle w(T(\beta)(\xi), Y\rangle = -\omega(T(1)(T(\beta)(\xi)), Y). 
$$
To prove the equality of these, we need 
to express right-translations in terms of the tangent groupoid multiplication. 

\begin{lem}
For $X\in A_m\Sigma$ and $g\in \Sigma$ with $\beta(g) = m$, 
\begin{equation}
\label{eq:0g}
X\sol \tilo_g = T(R_g)(X).  
\end{equation}
\end{lem}

\pf
Write $X$ as the derivative at $t = 0$ of a curve $h_t\in\Sigma_m$ with
$h_0 = 1_m$. One can likewise write $\tilo_g$ as the derivative of the curve constant 
at $g$. Then $X\sol \tilo_g$ is the derivative of the curve $h_tg$ and this 
is $T(R_g)(X)$. 

We now complete the proof that $\tilbeta\circ \omega^\flat = w\circ T(\beta)$. We
must show that
$$
\omega(\xi, T(R_g)(Y)) = \omega(T(1)(T(\beta)(\xi)), Y). 
$$
Using (\ref{eq:isot}) and (\ref{eq:0g}) the LHS is equal to 
$$
\omega(T(1)T(\beta)(\xi)\sol\xi, Y\sol \tilo_g)
= \omega(T(1)T(\beta)(\xi), Y) + \omega(\xi, \tilo_g)
$$
and $\omega(\xi, \tilo_g)$ is zero, since $\tilo_g$ is the zero of $T_g\Sigma$. 

The proof that $\tilalpha\circ \omega^\flat = w\circ T(\alpha)$ is similar. 
\pfend

Finally, the morphism property itself is (\ref{eq:wd}). Note that $w$ is an
isomorphism of vector bundles since $\omega^\flat$ is. 

\subsection*{Poisson structure on $P$}

Consider Figure \ref{fig:TGT*G}(a). This diagram and (b) are not diagrams of
morphisms (as in (\ref{eq:omegaflat})) but are single mathematical objects, part of 
the structure of which is indicated by the arrows. (See \S\ref{sect:da}.)

\begin{figure}[h]
\centering
\subfloat[]
{\xymatrix@=1.2cm{
T\Sigma \ar[r]^{\tilq}
  \ar@<-0.5ex>[d]\ar@<0.5ex>[d]
& \Sigma \ar@<-0.5ex>[d]\ar@<0.5ex>[d]\\
TP \ar[r]_p & P\\
}}
\qquad
\subfloat[]
{\xymatrix@=1.2cm{
T^*\Sigma \ar[r]^{\tilc}
  \ar@<-0.5ex>[d]\ar@<0.5ex>[d]
& \Sigma \ar@<-0.5ex>[d]\ar@<0.5ex>[d]\\
A^*\Sigma \ar[r]_q & P\\
}}
\caption{\label{fig:TGT*G}}
\end{figure}

Recall that the Lie algebroid of $\Sigma$, as a vector bundle, is the intersection of 
the kernel of
the bundle projection $\tilq$, which is the restriction of $T\Sigma$ to the identity
elements, and the kernel of $T(\alpha)$, which is the vector bundle tangent to the
source fibres. The bracket on $\Ga A\Sigma$ is given by associating to a section
$X$ of $A\Sigma$ a right-invariant vector field $\Ri{X}$, and showing that the
bracket of right-invariant vector fields is right-invariant. We now apply this
process to the cotangent structure in Figure \ref{fig:TGT*G}(b). 

First find the kernel of $\tilalpha$. To do this we use the following lemma. 

\begin{lem}
\label{lem:little}
Let $\eta\in T_h\Sigma$ and $\xi\in T_g\Sigma$ have $T(\alpha)(\eta) = 0 = 
T(\beta)(\xi)$. Then 
$$
\eta\sol\xi = T(L_h)(\xi) + T(R_g)(\eta).
$$
\end{lem}

\pf
In this case $\eta$ is tangent to $\alpha^{-1}(m)$ and $\xi$ is tangent to
$\beta^{-1}(m)$ where $m = \alpha h = \beta g$, and so $\eta\sol\xi$ can be
calculated from the restriction of multiplication to 
$\alpha^{-1}(m)\times \beta^{-1}(m)\to \Sigma$, which
does not involve a pullback. 
\pfend

Now take $X\in A\Sigma$. Then, applying the lemma, we have 
$$
X\sol (T(1)(aX) - X) = X + (T(1)(aX) - X) = T(1)(aX)
$$
and so 
$$
X^{-1} = T(1)(aX) - X.
$$
So we can rewrite the definition of $\tilalpha$ as
$$
\langle\tilalpha(\Phi), X\rangle = -\langle\Phi, T(L_g)(X^{-1})\rangle, 
$$
and the kernel of $\tilalpha$ is therefore the annihilator of $T^\beta\Sigma$. 

Suppose $T(\beta)(\xi) = 0$, where $\xi\in T_g\Sigma$, Then $T(\alpha)(\xi^{-1}) = 0$. 
So $\xi^{-1} = T(R_{g^{-1}})(X)$ for some $X\in A\Sigma$. That is, 
$\xi^{-1} = X\sol \tilo_{g^{-1}}$. Thence $\xi = \tilo_g\sol X^{-1}$.   

Clearly all pullbacks $\beta^*\mu$, for $\mu\in T^*P$, annihilate $T^\beta\Sigma$. 
Conversely, suppose that $\Phi\in T^*_{1_m}\Sigma$ annihilates $T^\beta\Sigma$. 
Then $\langle\Phi,X\rangle = \langle\Phi, T(1)(aX)\rangle$ for all $X\in A_m\Sigma$. Define
$\mu\in T^*_mP$ by $\langle\mu, x\rangle = \langle \Phi, T(1)(x)\rangle$ for
$x\in T_mP$. Then, for all $T(1)(x) + X\in T_{1_m}\Sigma$,  
$$ 
\langle\Phi, T(1)(x) + X\rangle = \langle \Phi, T(1)(x) + T(1)(aX)\rangle =
\langle\mu, x+aX\rangle = \langle \beta^*\mu, T(1)(x) + X\rangle, 
$$ 
so $\Phi = \beta^*\mu$. Thus the bundle for Figure \ref{fig:TGT*G}(b) which corresponds 
to $A\Sigma$ in (a), is $T^*P$. 

Now take $\mu\in\Om^1(P)$. We define $\Ri{\mu}\in\Om^1(\Sigma)$ by
$$
\Ri{\mu}(g) = \mu(\beta g)\sol\hato_g, 
$$
where $\hato_g$ is the zero element of $T^*_gG$. In fact $\mu(\beta g)\sol\hato_g$ is 
merely the pullback of $\mu(\beta g)$ to $g$. 

Take $\mu\in T^*_mP$ and $g$ with $\beta g = m$. Then
$$
\langle\mu\sol\hato_g, \eta\sol\xi\rangle = \langle\mu, \eta\rangle + 0
$$  
where $\eta = T(1)(y) + Y$ and $T(\beta)(\xi) = T(\alpha)(\eta) = y$. 
Now, regarding $\mu$ as in $T^*_{1_m}\Sigma$, we have
$$
\langle\mu, \eta\rangle = \langle\mu, T(\beta)(\eta)\rangle = \langle\mu,  y + aY\rangle.
$$
On the other hand, 
$$
\langle\beta^*_g\mu, \eta\sol\xi\rangle = \langle \mu, T(\beta)(\eta)\rangle 
= \langle\mu, y + aY\rangle.
$$

We can now define a bracket on $\Om^1(P)$ by
\begin{equation}
\label{eq:b1f}
\Ri{[\mu_1, \mu_2]} = [\Ri{\mu_1}, \Ri{\mu_2}],
\end{equation}
where the bracket on the right is the bracket on $\Om^1(\Sigma)$ transported
via $\omega^\flat$ from $\vf{\Sigma}$. The bracket on $\Om^1(P)$ is therefore 
skew-symmetric and satisfies the Jacobi identity. % thanks YKS
Note that (\ref{eq:b1f}) can also be
written as 
\begin{equation}
\label{eq:b1f-v2}
\beta^*[\mu_1, \mu_2] = [\beta^*\mu_1, \beta^*\mu_2],
\end{equation}

Regarding $A\Sigma\subseteq T\Sigma$ and $T^*P\subseteq T^*\Sigma$ as the
intersections of the relevant kernels, it follows from the fact that
$\omega^\flat$ is an isomorphism of vector bundles, that the restriction
to $A\Sigma\to T^*P$ is also. Denote this restriction temporarily
by $r$. Then, 
\begin{equation}
r = -w^*.   
\end{equation}
To see this, take $X\in A\Sigma$ and $y\in TP$. Then
$$
\langle w^*(X), y\rangle = \langle w(y), X\rangle =
- \omega(T(1)(y), X), 
$$
whereas 
$$
\langle r(X), y\rangle = \langle \omega^\flat(X), T(1)(y)\rangle 
= - \omega(X, T(1)(y)). 
$$
The brackets in $\Ga A\Sigma$ and $\Ga T^*P$ are defined in terms of
those on $\vf{\Sigma}$ and $\Om^1(\Sigma)$ by analogous processes and it
follows that $r$ is an isomorphism of Lie algebroids, to $T^*P$ with
the bracket defined by (\ref{eq:b1f}) and the anchor $a\circ r^{-1}$, 
where $a$ is the anchor of $A\Sigma$. 

\begin{lem}
$$
(a\circ r^{-1})^* = -a\circ r^{-1}.   
$$
\end{lem}

\pf
We first prove that $w\circ a  = a^*\circ r$. For $X,Y\in A\Sigma$, 
$$
\langle (w\circ a)(X),Y\rangle = \omega(T(1)(aX), Y)
$$
and 
$$
\langle (a^*\circ r)(X),Y\rangle 
= \langle r(X), a(Y)\rangle 
= - \langle w^*(X), a(Y)\rangle 
= - \langle w(aY), X\rangle 
= - \omega(T(1)(aY), X\rangle.
$$
We know that $T(1)(aX) = X+X^{-1}$ and $T(1)(aY) = Y + Y^{-1}$. So we have
$$
\langle (w\circ a)(X),Y\rangle = \omega(X, Y) + \omega(X^{-1},Y)\ \text{and}\ 
\langle (a^*\circ r)(X),Y\rangle = \omega(X,Y) + \omega(X,Y^{-1}).
$$
We must show that $\omega(X^{-1},Y) = -\omega(X,Y^{-1})= 0$. 

Take $X\in A\Sigma$. Since $r$ maps $X$ to $r(X)\in T^*P$, it follows that
$\omega^\flat$ maps $\Rarr{X}$ to $\Rarr{r(X)}$; that is, 
\begin{equation}
\label{eq:basic}
\omega^\flat(\Rarr{X}) = \beta^*(r(X)) = -\beta^*(w^*(X)).   
\end{equation}
So for any $\eta\in T\Sigma$, $\omega(\Rarr{X},\eta) = -\langle \beta^*(w^*(X)),\eta\rangle 
= -\langle w^*(X), T(\beta)(\eta)\rangle$. In particular if $T(\beta)(\eta) = 0$
then $\omega(\Rarr{X},\eta) = 0$. This completes the proof that 
$w\circ a  = a^*\circ r$.

Rewrite this as $a\circ r^{-1} = w^{-1}\circ a^*$. Now
$$
(a\circ r^{-1})^* = (w^{-1}\circ a^*)^* = 
a\circ (w^*)^{-1} = - a\circ r^{-1},
$$
as required. 
\pfend

Note from the proof that we have also shown that:

\begin{prop}
For $\xi,\eta\in T\Sigma$, if $T(\alpha)(\xi) = 0$ and $T(\beta)(\eta) = 0$ then
$\omega(\xi,\eta) = 0$.   
\end{prop}

Now define a Poisson structure $\pi$ on $P$ by 
\begin{equation}
\{f_1,f_2\} = \langle df_2, (a\circ r^{-1})(df_1)\rangle;  
\end{equation}
that is, $\pi^\# = a\circ r^{-1}$. That $\pi^\#$ is skew-symmetric has just
been proved, and the Jacobi identity follows from the isomorphism of $T^*P$ with
$A\Sigma$. So $\pi$ is a Poisson structure on $P$ and, by (\ref{eq:b1f-v2}), $\beta$
is a Poisson map. Since $\beta$ is a surjective submersion, it is the only such
Poisson structure. 

This completes the proof of (iii) and (iv) of Theorem \ref{thm:sgpd}. 

\begin{prop}
\label{prop:tangentP}
The diffeomorphism $-w\co TP\to A^*\Sigma$ is a Poisson map from the tangent lift 
Poisson structure on $TP$ to the Poisson structure dual to the Lie algebroid
structure on $A\Sigma$.  
\end{prop}

The \emph{tangent lift Poisson structure on $TP$} can be defined as the Poisson
structure which is dual to the Lie algebroid $T^*P$. It was originally defined
by T.~Courant \cite{Courant:1994} directly in terms of functions on $TP$. For $f\in\sfn{P}$
denote the function $TP\to\R$ corresponding to $df\co P\to T^*P$ by $\ell_{df}$. 
Then it is sufficient to define the Poisson brackets for all $\ell_{df}$ and
all $p^*f$, where $p\co TP\to P$ is the bundle projection. As in \cite{Courant:1994}, 
we define:
\begin{equation}
\{\ell_{df_1},\ell_{df_2}\} = \ell_{d\{f_1,f_2\}},\qquad
\{\ell_{df_1},p^*f_2\} = p^*\{f_1,f_2\},\qquad
\{p^*f_1,p^*f_2\} = 0. 
\end{equation}

Proposition \ref{prop:tangentP} now follows from the following general result. 
The proof is a simple exercise. 

\begin{prop}
Let $A$ and $B$ be Lie algebroids on the same base $M$, and let $\phi\co A\to B$
be a morphism of vector bundles over $M$. Then $\phi$ is a morphism of Lie algebroids
if and only if $\phi^*\co B^*\to A^*$ is a Poisson map. 
\end{prop}

\section{Midword} % 
\label{sect:midword}

We have given the proof that the base manifold of a symplectic groupoid has a
Poisson structure such that the target is a Poisson map, and that the
cotangent Lie algebroid of this Poisson structure is canonically
isomorphic to the Lie algebroid of the Lie groupoid. This is a very
striking result. It shows, in particular, that the Poisson structure
on the base manifold determines the symplectic groupoid up to local isomorphism. 

Theorem \ref{thm:sgpd} is the apotheosis of a classical question, 
as to whether a given Poisson manifold can be realized as the quotient 
of a symplectic manifold. In the example $T^*G\gpd\gog^*$ for $G$ a Lie
group, the Poisson manifold is the quotient of $T^*G$ under an
action of $G$. In the general case of Theorem \ref{thm:sgpd}, the symplectic
manifold $\Sigma$ is quotiented to its base $P$ by the surjective submersion
$\beta\co\Sigma\to P$. This may be regarded as the quotient over the action 
of the equivalence relation $\Sigma\times_\beta\Sigma$ on $\Sigma$. 

The classical question is the `realizability problem': given a Poisson manifold
$P$, is there a symplectic manifold $M$ and a surjective submersion 
$M\to P$ which is a Poisson map? In fact, under mild conditions, such a map can 
be modified to give a symplectic groupoid structure. 

The answer to the question --- \emph{given a Poisson manifold $P$ is there a symplectic
groupoid $\Sigma\gpd P$~?} --- was provided by Crainic and Fernandes 
\cite{CrainicF:2003,CrainicF:2004}, building on work of Cattaneo and Felder 
\cite{CattaneoF:2001}. 

The modern theory of symplectic realizations was introduced by Weinstein 
\cite{Weinstein:1987,CosteDW} in order to provide a route for 
the quantization of Poisson manifolds. 

In this article, we are concerned with the structures arising from symplectic 
groupoids. These are, I believe, most easily understood by studying Poisson
groupoids, of which symplectic groupoids are a particular case. Poisson groupoids
were introduced by Weinstein \cite{Weinstein:1988} and gave rise to the theory
of Lie bialgebroids. It is to some aspects of this theory that we now turn. 

\section{Poisson groupoids}
\label{sect:pgpds}

Let $P$ be a Poisson manifold. A closed submanifold $C\subseteq P$ is \emph{coisotropic}
if the Poisson anchor $\pi^\#\co T^*P\to TP,$ when
restricted to $(TC)^\circ,$ goes into $TC.$ Here $(TC)^\circ$ is the \emph{annihilator of
$TC$ in $T^*P$}; it is isomorphic to the conormal bundle $T_CP/TC.$

In terms of the bracket of functions, a closed submanifold $C$ is coisotropic in $P$
if, whenever $f,g\in\sfn{P}$ vanish on $C$, their bracket $\{f,g\}$ does also. 

Coisotropic submanifolds play the role that in symplectic geometry is played
by Lagrangian submanifolds. In particular, a map $\phi\co P\to Q$ of Poisson manifolds 
is a Poisson map if and only if the graph of $\phi$ is a coisotropic submanifold of 
$\Bar{Q}\times P$. 

\begin{df}[\cite{Weinstein:1988}]
\label{df:pgpd}
A \emph{Poisson groupoid} is a Lie groupoid $\sg\gpd P$ together with a Poisson 
structure on $\sg$ such that the graph of multiplication
\begin{equation}
\label{eq:graphP}
Gr = \{(hg,h,g)\st \alpha h = \beta g\}
\end{equation}
is coisotropic in $\Bar{\sg}\times\sg\times\sg$. 
\end{df}

Using (\ref{eq:tgraph}), $T(Gr)^\circ$ consists of triples $(\Theta,\Psi,\Phi)$
of elements of $T^*\sg$ such that
\begin{equation}
\label{eq:anncond}
\Theta(\eta\sol\xi) + \Psi(\eta) + \Phi(\xi) = 0
\end{equation}
for all $\eta, \xi\in T\sg$ such that $T(\alpha)(\eta) = T(\beta)(\xi)$. 
We claim that this implies $\Theta = -\Psi\sol\Phi$. First we must show that
$\tilalpha(\Psi) = \tilbeta(\Phi)$. Let $h,g\in\sg$ be the base elements of
$\eta, \xi$ and write $m = \beta(g)$. For $Y\in A_m\sg$ we have
$$
(\tilbeta(\Phi))(Y) = \Phi(T(R_g)(Y)) = \Phi(Y\sol\tilo_g).
$$
From (\ref{eq:anncond}) we now get
$$
\Phi(Y\sol\tilo_g) = \Theta(\tilo_h\sol\tilo_g) - \Psi(\tilo_h\sol Y^{-1}) 
= - \Psi(T(L_h)(Y^{-1})) 
$$
However $Y^{-1} = T(1)(aY) - Y$, since $Y\in A\sg$, and so we have
$$
(\tilbeta(\Phi))(Y) = \Psi(T(L_h)(Y-T(1)(aY))) = (\tilalpha(\Psi))(Y).  
$$
Now (\ref{eq:anncond}) can be written
\begin{equation}
\label{eq:anncond2}
T(Gr)^\circ = \{(-\Psi\sol\Phi,\Psi,\Phi)\st \tilalpha(\Psi) = \tilbeta(\Phi)\}. 
\end{equation}
The coisotropy condition now is that 
\begin{equation}
\label{eq:morphism}
\pi^\#(\Psi\sol\Phi) = \pi^\#(\Psi) \sol \pi^\#(\Phi)
\end{equation}
for $\Psi,\Phi\in T^*\sg$ with $\tilalpha(\Psi) = \tilbeta(\Phi)$. 
This is evidently a morphism condition on $\pi^\#$. Denote the base map
$A^*\sg\to TP$ temporarily by $b$. Then, for $\phi\in A^*\sg$, 
\begin{equation}
\label{eq:basemap}
\pi^\#(\tilone_\phi) = T(1)(b(\phi)). 
\end{equation}
From this it follows that $1_P$ is coisotropic in $\sg$. For if $\Phi\in T^*_{1_m}\sg$
annihilates $T_{1_m}(1_P)$ then, by (\ref{eq:cotid}), $\Phi = \tilone_\phi$ for some
$\phi\in A^*\sg$ and (\ref{eq:basemap}) is then the coisotropy condition for $1_P$. 

That $1_P$ is coisotropic in $\sg$ implies, as shown by Weinstein 
\cite{Weinstein:1988}, that $A^*\sg$ inherits a Lie algebroid structure from $T^*\sg$, 
the anchor of which is the restriction of $\pi^\#$; that is, $b$. We therefore write
$a_*$ for $b$ in what follows. 

To define the bracket we must verify two conditions, namely that 
\begin{itemize}
\item if $\Phi, \Psi\in\Ga T^*\sg$ have $\Phi|_{P},\ \Psi|_{P}\in\Ga A^*\sg$ then
$[\Phi, \Psi]|_{P}\in\Ga A^*\sg$ also;
\item if $\Phi,\Psi\in\Ga T^*\sg$ have $\Phi|_{P} = 0$ and 
$\Psi|_{P}\in\Ga A^*\sg,$ then $[\Phi,\Psi]|_{P} = 0.$
\end{itemize}

For the first, take any $\xi\in\vf{P}$ and extend it to a vector field on $\sg,$
also denoted $\xi.$ Take $\Phi, \Psi\in\Ga T^*\sg$ such that 
$\Phi|_P,\Psi|_P$ are sections of $(TP)^\circ.$ From (\ref{eq:ldf}) we have
\begin{equation}
\label{eq:davidt}
\langle[\Phi, \Psi], \xi\rangle
= \langle\ld_{\Phi^\#}(\Psi), \xi\rangle 
  -  \langle\ld_{\Psi^\#}(\Phi), \xi\rangle 
  - \langle d(\pi(\Phi, \Psi)), \xi\rangle.
\end{equation}
The first term is
\begin{equation}
\label{eq:davidta}
\langle\ld_{\Phi^\#}(\Psi), \xi\rangle
	= \ld_{\Phi^\#}(\langle\Psi, \xi\rangle) 
		- \langle{\Psi}, [{\Phi^\#}, \xi]\rangle.
\end{equation}
On the RHS, the first term is zero on $P$ because 
$\langle\Psi, \xi\rangle$ is constant on $P$ and $\Phi^\#$
is tangent to $P$. In the second term, both $\Phi^\#$ and $\xi$ 
are tangent to $P$, and so their 
Lie bracket is also; hence the pairing with $\Psi$ is zero on $P.$

The second term of (\ref{eq:davidt}) is likewise zero on $P.$ 
The third term is
$-\xi(\pi(\Phi, \Psi))$ and $\pi(\Phi, \Psi) 
= \langle \Psi, \Phi^\#\rangle$
vanishes on $P$ since $\Phi^\#$ is tangent to $P$.
Since $\xi\in\vf{P}$ was arbitrary, $[\Phi,\Psi]$ is a
section of $(TP)^\circ$.

To verify the second itemized condition, take $\Phi,\Psi\in\Ga T^*\sg$ with 
$\Phi|_P = 0$ and $\Psi|_P\in\Ga(TP)^\circ$. Let $\xi$ be any vector
field on $P$. 

Now in (\ref{eq:davidta}), the first term on the RHS is zero on $P$ 
because $\Phi^\#$ is zero on $P$. The second term,  
$- \langle\Psi, [\Phi^\#, \xi]\rangle$, may not be zero; see below. 

In the corresponding equation for $\langle\ld_{\Psi^\#}(\Phi), \xi\rangle$, 
the first term on the RHS is zero on $P$ because the bracket is zero on $P$
and $\Psi^\#$ is tangent to $P$, by coisotropy. The second term is 
zero on $P$ because $\Phi$ is zero on $P$. 

The third term is equal to $\ld_\xi(\langle\Psi, \Phi^\#\rangle) = 
\langle\ld_\xi(\Psi), \Phi^\#\rangle + \langle\Psi, [\xi, \Phi^\#]\rangle$.
The first term on the RHS is zero on $P$ because $\Phi^\#$ is so. The
second term may be nonzero, and cancels with the term above. 

Since $\xi$ was any vector field on $\sg$, this proves that $[\Phi,\Psi]$ is
zero on $P$. 

This shows that $(TP)^\circ$ is a Lie subalgebroid of $T^*\sg$. It is also
possible to regard the Lie algebroid structure on $A^*\sg$ as a quotient of
the Lie algebroid $T^*\sg$. 

We have now proved part of the following theorem, which should be compared 
with Theorem \ref{thm:sgpd}. 
\begin{thm}
\renewcommand{\theenumi}{\rm(\roman{enumi})}
\label{thm:pgpd}
Let $\sg\gpd P$ be a Poisson groupoid. Then:
\begin{enumerate}
\item The Poisson anchor $\pi^\#\co T^*\sg\to T\sg$ is a morphism of Lie
groupoids with base map denoted $a_*\co A^*\sg\to TP$. 
\item $1_P$ is coisotropic in $\sg$.
\item The Lie algebroid structure on $T^*\sg$ induces a Lie algebroid structure on the
dual bundle $A^*\sg$ for which the anchor is $a_*$. 
\item Inversion is an antiPoisson map. 
\item There is a unique Poisson structure on $P$ such that $\beta$ is a Poisson map.
\end{enumerate}
\end{thm}

The chief difference between Theorem \ref{thm:pgpd} and Theorem \ref{thm:sgpd} 
concerns the Lie algebroid
structure on $A^*\sg$. In \ref{thm:sgpd}(iv) there is a natural isomorphism 
of Lie algebroids $A\Sigma\cong T^*P$. The dual of this, $TP\to A^*\Sigma$, 
may be used to put a Lie algebroid structure on $A^*\Sigma$, though there is
little reason to do this if one is working only with symplectic groupoids. 
  
The proof of Theorem \ref{thm:pgpd}(iv) is similar to that for symplectic groupoids. 
Note however that whereas inversion $\sg\to\sg$ is an antiPoisson diffeomorphism, the 
groupoid inversion
$T^*\sg\to T^*\sg$ is a Lie algebroid isomorphism, without change of sign. This confirms,
in a small way, the naturality of working with the cotangent structures. 

The Poisson structure on $P$ is defined by $\pi_P^\# = a_*\circ a^* = - a\circ a^*_*$. 
It follows that $\beta\co \sg\to P$ is a Poisson map.  

We need to explicate further the relationship between the Lie algebroid structures on
$A\sg$ and $A^*\sg$; in particular we want to express the relationship between $A\sg$
and $A^*\sg$ without reference to the underlying Poisson groupoid. 

With any Lie algebroid $A$ on base manifold $M$ there is a cochain complex 
$$
\xymatrix@=8mm{
\sfn{M}\ar[r]^d & \Ga (A^*)\ar[r]^d & \Ga\extt{2}{A^*}\ar[r]^d & 
\Ga\extt{3}{A^*} \ar[r]^d & \cdots\\
}
$$
defined by the natural extension of the coboundary operator for de Rham
cohomology and for Lie algebra cohomology. 

\begin{thm}
\label{thm:lbia}
Let $\sg\gpd P$ be a Poisson groupoid. Then for the coboundary operator $d_*$ of the   
Lie algebroid $A^*\sg$ defined above we have, for $X,Y\in\Ga A\sg$, 
\begin{equation}
\label{eq:lbia}
d_*[X,Y] = [X, d_*Y] + [d_*X,Y].   
\end{equation}
\end{thm}

We will not give the proof here; see \cite{MackenzieX:1994} 
or \cite[\S12.1]{Mackenzie:GT}. 

On the RHS we have the bracket of sections of $A\sg$ with sections of 
$\Ga\extt{2}{A\sg}$. This bracket is defined in terms of decomposable
elements by 
$$
[X,Y\wedge Z] = [X,Y]\wedge Z + Y\wedge[X,Z]. 
$$
and $[X,\eta] = -[\eta,X]$ for $X\in\Ga A\sg$ and $\eta\in\Ga\extt{2}{A\sg}$. 
The bracket on $\Ga A\sg$ may be extended to the exterior algebra $\Ga\extt{}{A\sg}$;
this is the \emph{Schouten} or \emph{Gerstenhaber bracket for} $A\sg$. Details
and references are given in \cite[\S7.5]{Mackenzie:GT}. 

\begin{df}
\label{df:lbia}
Let $A$ be a Lie algebroid on $M$ and suppose that $A^*$ has a Lie algebroid
structure. Then $(A, A^*)$ is a \emph{Lie bialgebroid} if (\ref{eq:lbia}) holds
for all $X, Y\in\Ga A$.   
\end{df}

Lie bialgebroids were defined by Ping Xu and the author in \cite{MackenzieX:1994}, 
following the construction
of the Lie algebroid structure on the dual by Weinstein \cite{Weinstein:1988}. 
Equation (\ref{eq:lbia}) was immediately shown by Kosmann-Schwarzbach 
\cite{Kosmann-Schwarzbach:1995} to extend to elements of arbitrary degree of 
$\Ga\extt{}{A}$; therefore $d_*$ is a derivation of the Schouten
bracket on $\Ga\extt{}{A}$. Equation (\ref{eq:lbia}) implies its dual form, 
$d[\phi,\psi]_* = [\phi, d\psi]_* + [d\phi,\psi]_*$, where $\phi,\psi\in\Ga A^*$
and $[~,~]_*$ is the bracket on $\Ga A^*$. 

Equation (\ref{eq:lbia}) is thus in many respects well-understood, but it remains
a very nonlinear equation. We will now show how the relationship between the Lie 
algebroid structures on $A\sg$ and $A^*\sg$ may be understood in a `diagrammatic' way. 

\section{The diagrammatic approach}
\label{sect:da}

We work with a given Poisson groupoid $\sg\gpd P$. The key to the diagrammatic 
approach is to work with the cotangent bundle $T^*\sg$ rather than $\sg$ itself. 
Since $\sg$ is a Poisson manifold, $T^*\sg$ is a Lie algebroid, and since $\sg$
is a Lie groupoid, $T^*\sg$ has a Lie groupoid structure. These structures are
shown in Figure~\ref{fig:VBgpds}(a).  

\begin{figure}[htb]
\begin{center}
\subfloat[]{
\xymatrix@=12mm{
T^*\sg  \ar[r]\ar@<-0.5ex>[d]\ar@<0.5ex>[d] & \sg\ar@<-0.5ex>[d]\ar@<0.5ex>[d]\\
A^*\sg \ar[r] & P\\
}}
\hspace*{30mm}
\subfloat[]{
\xymatrix@=12mm{
\Om  \ar[r]\ar@<-0.5ex>[d]\ar@<0.5ex>[d] & \sg\ar@<-0.5ex>[d]\ar@<0.5ex>[d]\\
A \ar[r] & M\\
}}
\hspace*{30mm}
\subfloat[]{
\xymatrix@=12mm{
T\sg  \ar[r]\ar@<-0.5ex>[d]\ar@<0.5ex>[d] & \sg\ar@<-0.5ex>[d]\ar@<0.5ex>[d]\\
TP \ar[r] & P\\
}}
\end{center}
\caption{\ \label{fig:VBgpds}}
\end{figure}

For any Lie groupoid $\sg\gpd P$, it is straightforward to show that $T^*\sg$ is a 
\VBgpd in the following sense.  

\begin{df}[\cite{Pradines:1988}]
\label{df:VBgpd}
A {\em \VBgpd}\ $(\Om;\sg,A;M)$ is a structure as shown in 
Figure~\ref{fig:VBgpds}(b) in which $\Om$ is a vector bundle 
over $\sg,$ which is a Lie groupoid over $M,$ and $\Om$ is also 
a Lie groupoid over $A,$ which is a vector bundle over $M,$ 
subject to the condition that the structure maps of the groupoid 
structure (source, target, identity, multiplication, inversion) 
are vector bundle morphisms, and the \lq double source map\rq\
$\Om\to \sg\times_{M}A$ formed from the bundle projection and the
source on $\Om$, is a surjective submersion.
\end{df}

According to the duality for \VBgpds introduced by Pradines \cite{Pradines:1988},
the \VBgpd $T^*\sg$ is dual to the \VBgpd $T\sg$ shown in Figure~\ref{fig:VBgpds}(c),
the tangent prolongation of $\sg\gpd P$. 

Each diagram in Figure~\ref{fig:VBgpds} denotes the full set of
structures just described, and constitutes a single mathematical 
object; they should not be read, for example, as showing a 
morphism of groupoids. 

For $\sg\gpd P$ a Poisson groupoid, each horizontal structure in 
Figure~\ref{fig:VBgpds}(a) is a Lie algebroid. These are compatible
with the vertical groupoid structures in the following sense. 

\begin{df}[\cite{Mackenzie:1992}]
\label{df:LAgpd}
An {\em \LAgpd}\ $(\Om;\sg,A;M)$ is a structure as shown in 
Figure~\ref{fig:VBgpds}(b) in which $\Om$ is a Lie algebroid
over $\sg,$ which is a Lie groupoid over $M,$ and $\Om$ is also 
a Lie groupoid over $A,$ which is a Lie algebroid over $M,$ 
subject to the condition that the structure maps of the groupoid 
structure (source, target, identity, multiplication, inversion) 
are Lie algebroid morphisms, and the \lq double source map\rq\
$\Om\to \sg\times_{M}A$ formed from the bundle projection and the
source on $\Om$, is a surjective submersion.
\end{df}

We proved in \S\ref{sect:pgpds} that the anchor $\pi^\#\co T^*\sg\to T\sg$
is a morphism of Lie groupoids; that result may equally be formulated as
stating that each of the groupoid structure maps preserves the anchors. 
It is necessary to also prove that the groupoid structure maps preserve
the brackets; for this proof, see \cite{Mackenzie:1992}. Thus 
$T^*\sg$, for $\sg\gpd P$ a Poisson groupoid, is an \LAgpd. This has
important ramifications: one may differentiate an \LAgpd, and one may
seek to integrate it. 

Consider first the result of differentiating the \LAgpd in Figure~\ref{fig:VBgpds}(b). 
Applying the Lie functor to the vertical groupoid structures we obtain a structure
as shown in Figure~\ref{fig:doubles}(a). Here each arrow represents a vector bundle 
structure; altogether these make $A\Om$ a double vector bundle, as defined by Pradines 
\cite{Pradines:DVB,Pradines:1974b}. 

\begin{figure}[htb]
\begin{center}
\subfloat[]{
\xymatrix@=12mm{
A\Om  \ar[r]\ar[d] & A\sg\ar[d]\\
A \ar[r] & M\\
}}\hfil
\subfloat[]
{\xymatrix@=12mm{
AT^*\sg \ar[d]\ar[r] & A\sg \ar[d]\\
A^*\sg \ar[r] & P\\
}}\hfil
\subfloat[]{
\xymatrix@=12mm{
D  \ar[r]\ar[d] & B\ar[d]\\
A \ar[r] & M\\
}}\hfil
\subfloat[]{
\xymatrix@=12mm{
T^*A  \ar[r]\ar[d] & A\ar[d]\\
A^* \ar[r] & M\\
}}
\end{center}
\caption{\ \label{fig:doubles}}
\end{figure}

The vertical structures are evidently Lie algebroids. It is also true, though not
so very quick to prove \cite{Mackenzie:2000}, that the horizontal structures inherit 
Lie algebroid structures by a process which is an extension of the construction of 
the tangent prolongation of a Lie algebroid. Thus all four sides of 
Figure~\ref{fig:doubles}(a) have Lie algebroid structures. 

Applying this to the cotangent \LAgpd of a Poisson groupoid, Figure~\ref{fig:VBgpds}(a), 
we obtain Figure~\ref{fig:doubles}(b). Here there is a canonical isomorphism
$AT^*\sg \cong T^*A\sg$ which preserves the Lie algebroid structures
\cite{MackenzieX:1994}. Thus the Lie algebroid structure on $AT^*\sg\to A\sg$ may be 
obtained as the cotangent Lie algebroid structure of the Poisson structure on $A\sg$ 
which is dual to the Lie algebroid structure on $A^*\sg$. 

When $\sg\gpd P$ is a Poisson Lie group $G$ (with $P = \{\cdot\}$), $T^*AG$ reduces as
a manifold to $\gog\times\gog^*$. The vertical Lie algebroid structure arises  
by lifting the Lie algebra structure on $\gog$ to the pullback vector bundle 
$\gog^*\times\gog\to\gog^*$; similarly with $\gog^*$. In turn the classical 
Drinfel'd double structure on $\gog\oplus\gog^*$ may be obtained from $T^*A\sg$ by
a process of diagonalization \cite{Mackenzie:2011}. Thus Figure~\ref{fig:doubles}(b)
may be regarded as an extension to Poisson groupoids of the classical Drinfel'd double
of a Poisson Lie group. 

The relationship between the Lie algebroid structures in diagrams 
such as those in Figure~\ref{fig:doubles}(a)(b) cannot be characterized in 
the same way as for \VBgpds and \LAgpds; it is not possible to say that the Lie 
algebroid brackets on the (say) vertical structures are morphisms with respect
to the horizontal structures. This problem can be solved by using the duality
theory of double vector bundles. Very briefly, given a double vector bundle
as in Figure~\ref{fig:doubles}(c) with Lie algebroid structures on each side,
assume that each Lie algebroid structure on $D$ is linear with respect to the 
other; then dualizing each structure on $D$ leads to a pair of Poisson vector
bundles in duality, and $D$ is defined to be a \emph{double Lie algebroid} if
the corresponding dual Lie algebroids form a Lie bialgebroid. A full account
is in \cite{Mackenzie:2011}. The case in which one pair of parallel Lie
algebroids is zero is variously called an \LAvb or a \VBalgd. Thus one may
regard a double Lie algebroid as a double vector bundle admitting a horizontal
\VBalgd structure and a vertical \VBalgd structure which are suitably compatible
\cite{Gracia-SazMehta:2010}. 

Any vector bundle $A\to M$ gives rise to a \emph{cotangent double vector bundle} 
$T^*A$ as shown in Figure~\ref{fig:doubles}(d). The structure $T^*A\to A^*$ arises
from using the canonical diffeomorphism $T^*A\cong T^*A^*$ \cite{MackenzieX:1994}. 
Now suppose that both $A$ and $A^*$ have Lie algebroid structures, not necessarily
related. The Lie algebroid structure on $A^*$ induces a Poisson structure on $A$ and 
this induces a Lie algebroid structure on $T^*A\to A$. Likewise the Lie algebroid 
structure on $A$ induces a Lie algebroid structure on $T^*A^*\to A^*$. These structures
make $T^*A\cong T^*A^*$ a double Lie algebroid if and only if $A$ and $A^*$ form a 
Lie bialgebroid \cite{Mackenzie:2011}. 

Thus it is reasonable to argue that the cotangent double Figure~\ref{fig:doubles}(d)
plays the role for Lie bialgebroids which is played for Lie bialgebras by the classical
Drinfel'd double. In particular the bialgebroid equation (\ref{eq:lbia}) is encapsulated 
in the double Lie algebroid conditions for the cotangent double $T^*A$. 

The notion of double Lie algebroid provides an alternative to the extension of the 
classical Drinfel'd double to Lie bialgebroids of Liu, Weinstein and Xu 
\cite{LiuWX:1997}. The theory of Courant algebroids has led to the important
concept of Dirac structures, whereas double Lie algebroids arise as second-order 
invariants of double Lie groupoids. We will describe this briefly below. 

Double Lie algebroids have been defined in terms of supergeometry by Th.~Voronov
\cite{Voronov:2012-QM}, and this is an exceptionally clear and elegant formulation. 
In the late 1990s, Va{\u\i}ntrob showed that a homological vector field of degree
$+1$ on a supermanifold is equivalent to a Lie algebroid structure on the parity
reversed bundle. A double Lie algebroid is then defined by two homological
vector fields, of suitable weights, on a double parity reversion, provided that
the homological vector fields \emph{commute}. This includes the construction of the super
cotangent double by Roytenberg \cite{Roytenberg:thesis}. 

Very recently a formulation of double Lie algebroids in terms of representations
up to homotopy has been given by Gracia-Saz, Jotz Lean, Mehta and the author 
\cite{GSJLMM}. 

It was shown in \cite{MackenzieX:2000} that given a Lie bialgebroid $(A, A^*)$
in which $A = A\sg$ is the Lie algebroid of a source-simply-connected Lie
groupoid $\sg$, there is a Poisson structure on $\sg$ making it a Poisson 
groupoid, and inducing the given Lie algebroid structure on $A^*$. In terms 
of double structures, this shows that the cotangent double Lie algebroid
$T^*A\sg$ may be integrated to an \LAgpd $T^*\sg$. Very recently Bursztyn, 
Cabrera and del Hoyo \cite{BursztynCDH} have given a general integrability 
result for double Lie algebroids in terms of \LAgpds. 

We turn now to the question of integrating an \LAgpd, such as that in 
Figure~\ref{fig:VBgpds}(a). It was shown by Lu and Weinstein \cite{LuW:1989}
that underlying any Poisson Lie group is a symplectic double groupoid. 

In general a \emph{double Lie groupoid} consists of a manifold $S$ with two
Lie groupoid structures, $S\gpd H$ and $S\gpd V$, where $H$ and $V$ are
themselves Lie groupoids over a manifold $M$, as shown in Figure~\ref{fig:dgpds}(a), 
such that the structure maps of each groupoid structure are morphisms 
with respect to the other structure, and satisfying a double source condition 
analogous to those in Definitions \ref{df:VBgpd} and \ref{df:LAgpd}. 

In the same way that elements of a groupoid are visualized as arrows, elements
of a double groupoid are visualized as squares, 
the horizontal edges of which come from $H$, the vertical edges from $V$, 
and the four corners from $M$; see Figure \ref{fig:dgpds}(b)(c). 
If $v'_1 = v_2$ one can compose the two elements shown horizontally --- 
the vertical edges of the composite, shown in Figure \ref{fig:dgpds}(d), 
are determined by the groupoid axioms for $S\gpd V$ and the horizontal edges by the
condition that the source and target for $S\gpd H$ are morphisms. 

The condition that each groupoid multiplication is a morphism with respect to the
other structure is an interchange law: given four squares which can be arranged 
to form a larger square with matching inner edges, the result of composing
vertically and then horizontally is the same as composing horizontally and then 
vertically. 

\begin{figure}[h]
\centering
\subfloat[]
{\xymatrix@=1.5cm{
S \ar@<0.5ex>[r]\ar@<-0.5ex>[r]  \ar@<-0.5ex>[d]\ar@<0.5ex>[d]
& V \ar@<-0.5ex>[d]\ar@<0.5ex>[d]\\
H \ar@<-0.5ex>[r]\ar@<0.5ex>[r] & M\\
}}\hfil
\subfloat[]
{\xymatrix@=7mm{ 
&&  \ar[ll]_{h'_2}  \\
&  & \\
 \ar[uu]^{v'_2} && \ar[ll]^{\rule[-2mm]{0pt}{5mm}h'_1} \ar[uu]_{v'_1} \\
}}\hfil
\subfloat[]
{\xymatrix@=7mm{ 
&&  \ar[ll]_{h_2}  \\
&  & \\
 \ar[uu]^{v_2} && \ar[ll]^{\rule[-2mm]{0pt}{5mm}h_1} \ar[uu]_{v_1} \\
}}\hfil
\subfloat[]
{\xymatrix@=7mm{ 
&&  \ar[ll]_{h'_2h_2}  \\
&  & \\
 \ar[uu]^{v'_2} && \ar[ll]^{\rule[-2mm]{0pt}{5mm}h'_1h_1} \ar[uu]_{v_1} \\
}}\hfil
\caption{\ \label{fig:dgpds}}
\end{figure}

Given a double Lie groupoid one may apply the Lie functor to the horizontal
structures and obtain an \LAgpd as in Figure~\ref{fig:firstdiff}(a). As described
earlier, the Lie functor may then be applied to the vertical structures to
give the double Lie algebroid $A_VA_HS$ of $S$; see Figure~\ref{fig:firstdiff}(b).
Equally, one may apply the Lie functor to the vertical structures and obtain an 
\LAgpd as in Figure~\ref{fig:firstdiff}(c), and then to the horizontal structures
to give the double Lie algebroid $A_HA_VS$ of $S$; see Figure~\ref{fig:firstdiff}(d).
There is a canonical diffeomorphism from $A_VA_HS$ to $A_HA_VS$. 
  
For example if $S = M^4$, one obtains the \LAgpd $TM\times TM$ and the double Lie 
algebroid is the iterated tangent bundle $T^2M$. 

\begin{figure}[h]
\centering
\subfloat[]
{\xymatrix@=1.5cm{
A_HS \ar@<0.5ex>[d]\ar@<-0.5ex>[d]\ar[r]
& V \ar@<0.5ex>[d]\ar@<-0.5ex>[d]\\
AH \ar[r] & M\\
}}\hfil
\subfloat[]
{\xymatrix@=15mm{
A_VA_HS \ar[r] \ar[d] & AV \ar[d] \\
AH \ar[r] & M 
}}\hfil
\subfloat[]
{\xymatrix@=1.5cm{
A_VS \ar@<0.5ex>[r]\ar@<-0.5ex>[r]  \ar[d]
& AV \ar[d]\\
H \ar@<-0.5ex>[r]\ar@<0.5ex>[r] & M\\
}}\hfil
\subfloat[]
{\xymatrix@=15mm{
A_HA_VS \ar[r] \ar[d] & AV \ar[d] \\
AH \ar[r] & M 
}}\hfil
\caption{\ \label{fig:firstdiff}}
\end{figure}

Now consider a Poisson Lie group $G$ and let $\Dd$ denote the simply-connected 
Lie group corresponding to the classical Drinfel'd double $\gog\bowtie\gog^*$. 
Then the inclusions $\gog\to\gog\bowtie\gog^*$ and $\gog^*\to\gog\bowtie\gog^*$ 
induce morphisms $G\to\Dd$, $g\mapsto\Bar{g}$ and $G^*\to\Dd$, $\phi\mapsto\Bar{\phi}$, 
where $G^*$ is the simply-connected Lie group corresponding to $\gog^*$. 
Lu and Weinstein \cite{LuW:1989} define a double Lie groupoid $S$ as shown in
Figure~\ref{fig:LuW}(a), for which the elements are quadruples $(g_2,g_1,\phi_2,\phi_1)
\in G\times G\times G^*\times G^*$ such that $\Bar{g_2}\,\Bar{\phi_1} = 
\Bar{\phi_2}\,\Bar{g_1}\in\Dd$, 
as shown in Figure~\ref{fig:LuW}(b). 

\begin{figure}[h]
\centering
\subfloat[]
{\xymatrix@=1.5cm{
S \ar@<0.5ex>[r]\ar@<-0.5ex>[r]  \ar@<-0.5ex>[d]\ar@<0.5ex>[d]
& G \ar@<-0.5ex>[d]\ar@<0.5ex>[d]\\
G^* \ar@<-0.5ex>[r]\ar@<0.5ex>[r] & \{\cdot\}\\
}}\hfil
\subfloat[]
{\xymatrix@=7mm{ 
&&  \ar[ll]_{\phi_2}  \\
&  & \\
 \ar[uu]^{g_2} && \ar[ll]^{\rule[-2mm]{0pt}{5mm}\phi_1} \ar[uu]_{g_1} \\
}}
\hfil
\subfloat[]
{\xymatrix@=1.5cm{
A_HS \ar@<0.5ex>[d]\ar@<-0.5ex>[d]\ar[r]
& G \ar@<0.5ex>[d]\ar@<-0.5ex>[d]\\
\gog^* \ar[r] & \{\cdot\}\\
}}\hfil
\caption{\ \label{fig:LuW}}
\end{figure}

The Lie algebra structure on $\gog\bowtie\gog^*$ induces a symplectic 
structure on $S$ which makes $S$ a symplectic groupoid with
respect to both groupoid structures. It then follows that $A_HS\cong T^*G$ 
as Lie algebroids and as Lie groupoids. In effect, \cite{LuW:1989} has 
integrated the \LAgpd in Figure~\ref{fig:VBgpds}(a), in the case of 
Poisson Lie groups, to the symplectic double groupoid $S$. 

A general integration result, starting with an \LAgpd $\Om$ as in 
Figure~\ref{fig:VBgpds}(b), where $\Om\to\sg$ is the Lie algebroid of a 
suitable Lie groupoid $\Ga\gpd \sg$, and which constructs a second Lie 
groupoid structure $\Ga\gpd\scrh$, with $A\scrh\cong A$, so that $\Ga$ 
becomes a double Lie groupoid over $\sg$ and $\scrh$, seems very far 
from being accessible. 

\section*{Endword}
\addcontentsline{toc}{section}{Endword}

In addition to works already cited, other valuable sources on 
Poisson geometry include the books by Vaisman \cite{Vaisman:LGPM}, 
Cannas da Silva and Weinstein \cite{CannasdaSilvaW:GMNA}, 
Dufour and Zung \cite{DufourZ:PSNF}, and 
Laurent-Gengoux et al \cite{Laurent-GengouxPV}. 

\addcontentsline{toc}{section}{References}

\newcommand{\noopsort}[1]{} \newcommand{\singleletter}[1]{#1} \def\cprime{$'$}
  \def\cprime{$'$}

\end{document}